\documentclass{amsart}

\newtheorem{theorem}{Theorem}[section]
\newtheorem{corollary}[theorem]{Corollary}
\newtheorem{lemma}[theorem]{Lemma}

\numberwithin{equation}{section}

\begin{document}

\title{Structure of Group Invariants \\ of a Quasiperiodic Flow}


\author{Lennard F. Bakker}

\email{\rm bakkerl@member.ams.org}


\begin{abstract}

The multiplier representation of the generalized symmetry group of a quasiperiodic flow on the n-torus defines, for each subgroup of the multiplier group of the flow, a group invariant of the smooth conjugacy class of that flow. This group invariant is the internal semidirect product of a subgroup isomorphic to the n-torus by a subgroup isomorphic to that subgroup of the multiplier group. Each subgroup of the multiplier group is a multiplicative group of algebraic integers of degree at most n, which group is isomorphic to an abelian group of n by n unimodular matrices. 

\end{abstract}

\maketitle

\markboth{\sc L.F. Bakker}{\sc Structure of Group Invariants}

\section{Introduction}

The generalized symmetry group, $S_\phi$, of a smooth (i.e. $C^\infty$) quasiperiodic flow $\phi:{\mathbb R}\times T^n\to T^n$, $n\geq 2$, is group theoretic normalizer of the abelian group of diffeomorphisms generated by $\phi$:
\[S_\phi=N_{\hbox{\rm Diff($T^n$)}}\big(F_\phi\big),\]
where $\hbox{\rm Diff}(T^n)$ is the group of diffeomorphisms of $T^n$, and
$F_\phi=\{\phi_t:t\in{\mathbb R}\}$. The quasiperiodic flow $\phi$ is generated by the vector field
\[X(\theta)=\frac{d}{dt}\phi_t(\theta)\bigg\vert_{t=0}.\]
An $R\in\hbox{\rm Diff}(T^n)$ which belongs to $S_\phi$ is characterized by a number in ${\mathbb R}^*={\mathbb R}\setminus\{0\}$. (Here and elsewhere, ${\bf T}$ is the tangent functor, and $R_*X={\bf T}RXR^{-1}$.)

\begin{theorem}\label{char}
The following are equivalent:

{\rm a)} $R\in S_\phi${\rm ;}

{\rm b)} there exists a unique $\alpha\in{\mathbb R}^*$ such that 
$R\phi_t=\phi_{\alpha t}R$ for all $t\in{\mathbb R}${\rm ;}

{\rm c)} there exists a unique $\alpha\in{\mathbb R}^*$ such that $R_*X=\alpha X$.
\end{theorem}

\begin{proof}
See Lemma 10.3 and Theorem 13.1 in \cite{BC}.
\end{proof}

The multiplier representation $\rho_\phi:S_\phi\to{\mathbb R}^*=\hbox{\rm GL}({\mathbb R})$ is a linear representation of $S_\phi$ in ${\mathbb R}$ (see \cite{BA2}) which takes an $R$ in $S_\phi$ to the unique number $\alpha$ appearing in parts b) and c) of Theorem \ref{char}. The image, $\rho_\phi(S_\phi)$, is the multiplier group of $\phi$, and is a subgroup of the abelian group $\hbox{\rm GL}({\mathbb R})$. For each $\Lambda<\rho_\phi(S_\phi)$, the multiplier representation induces the short exact sequence of groups,
\[ \hbox{\rm id}_{T^n}\to\hbox{\rm ker\ }\rho_\phi\to \rho_\phi^{-1}(\Lambda)\to \Lambda\to 1,\]
in which $\hbox{\rm id}_{T^n}$ is the identity diffeomorphism of $T^n$,
$\hbox{\rm ker}\rho_\phi\to\rho_\phi^{-1}(\Lambda)$ is the canonical monomorphism, and $j_\Lambda:\rho_\phi^{-1}(\Lambda)\to\Lambda$ is the restriction of $\rho_\phi$ to $\rho_\phi^{-1}(\Lambda)$. (By the Fundamental Theorem on Homomorphisms (p. 10 \cite{AB}), $\Lambda\cong\rho_\phi^{-1}(\Lambda)/\hbox{\rm ker\ }\rho_\phi$.) This states that $\rho_\phi^{-1}(\Lambda)$ is a group extension of $\hbox{\rm ker}\rho_\phi$ by the abelian group $\Lambda$. It will be shown that 
this extension splits for any $\Lambda<\rho_\phi(S_\phi)$, that $\hbox{\rm ker}\rho_\phi\cong T^n$, and that $\Lambda$ is a multiplicative group of real algebraic integers of degree at most $n$, which is isomorphic to an abelian subgroup of $\hbox{\rm GL}(n,{\mathbb Z})$.

\section{Multipliers and Quasiperiodic Flows}

A flow $\phi$ on $T^n$ is quasiperiodic if and only if there exists a $V\in\hbox{\rm Diff}(T^n)$ such that $Y=V_*X$ is a constant vector field whose coefficients are independent over ${\mathbb Q}$ (see pp. 79-80 \cite{HB}). Real numbers $a_1,a_2,...,a_n$ are independent over ${\mathbb Q}$ if for $m=(m_1,m_2,...,m_n)\in{\mathbb Z}^n$, the equation
\[\sum^n_{j=1}m_ja_j=0\]
implies that $m_j=0$ for all $j=1,2,...,n$.

The nonzero number $\rho_\phi(R)$ is called the multiplier of $R\in S_\phi$. An $R\in S_\phi$ with $\rho_\phi(R)=1$ is known as a (classical) symmetry of $\phi$ (p. 1 \cite{AR}); the symmetry group of $\phi$ is $\hbox{\rm ker\ }\rho_\phi=\rho_\phi^{-1}(\{1\})$. An $R\in S_\phi$ with $\rho_\phi(R)=-1$ is called a reversing symmetry (p. 4 \cite{LR}); if $R^2=\hbox{\rm id}_{T^n}$, then $R$ is a reversing involution or a classical time-reversing symmetry of $\phi$; the reversing symmetry group of $\phi$ is $\rho_\phi^{-1}(\{1,-1\})$ (p. 8 \cite{LR}). An $R\in S_\phi$ with $\rho_\phi(R)\ne\pm1$, if it exists, is a new type of symmetry of $\phi$, called a generalized symmetry of $\phi$. Generalized symmetries are known to exist for quasiperiodic flows whose frequencies satisfy certain algebraic relationships (see \cite{BA1} for those algebraic relationships and examples on $T^2$ and $T^3$).

\begin{theorem}\label{mult}
If $\phi$ is a quasiperiodic, then $\{1,-1\}<\rho_\phi(S_\phi)$.
\end{theorem}

\begin{proof}
Suppose $\phi$ is quasiperiodic. Then there is a $V\in\hbox{\rm Diff}(T^n)$ such that $Y=V_*X$ is a constant vector field. Let $\psi$ be the flow generated by $Y$. For any $t\in{\mathbb R}$, the diffeomorphism $\psi_t$ satisfies $(\psi_t)_*Y=Y$, so that $1\in\rho_\psi(S_\psi)$. On the other hand, the map $N:T^n\to T^n$ defined by $N(\theta)=-\theta$ satisfies $N_*Y=-Y$, so that $-1\in\rho_\psi(S_\psi)$. Because $Y=V_*X$, the flows $\phi$ and $\psi$ are smoothly conjugate. This implies that $\rho_\phi(S_\phi)=\rho_\psi(S_\psi)$ (Theorem 12.2 \cite{BC}), and so $\{1,-1\}<\rho_\phi(S_\phi)$.
\end{proof}

\begin{theorem}\label{nonabelian}
If $\phi$ is quasiperiodic and $\{1\}\ne\Lambda<\rho_\phi(S_\phi)$, then $\rho_\phi^{-1}(\Lambda)$ is nonabelian, and hence the generalized symmetry group of $\phi$ and the reversing symmetry group of $\phi$ are nonabelian.
\end{theorem}

\begin{proof}
Suppose $\phi$ is quasiperiodic and $\{1\}\ne\Lambda$. Then there is an $R\in S_\phi$ such that $\alpha=\rho_\phi(R)\ne1$. By Theorem \ref{char}, $R\phi_1=\phi_{\alpha}R$. If $\phi_1=\phi_{\alpha}$, then $\phi$ would be periodic. Thus, $\rho_\phi^{-1}(\Lambda)$ is nonabelian. By Theorem \ref{mult}, both $\rho_\phi(S_\phi)$ and $\rho_\phi\big(\rho_\phi^{-1}(\{1,-1\})\big)$ contain $-1$, so that $S_\phi=\rho_\phi^{-1}(\rho_\phi(S_\phi))$ and $\rho_\phi^{-1}(\{1,-1\})$ are both nonabelian.
\end{proof}

For any $\Lambda<\rho_\phi(S_\phi)$, $\rho_\phi^{-1}(\Lambda)$ is an invariant of the smooth conjugacy class of $\phi$ in the sense that if $\phi$ and $\psi$ are smoothly conjugate, then $\rho_\phi^{-1}(\Lambda)$ and $\rho_\psi^{-1}(\Lambda)$ are conjugate subgroups of $\hbox{\rm Diff}(T^n)$ (Theorem 13.3 \cite{BC}). Because a quasiperiodic flow $\phi$ is smoothly conjugate to a quasiperiodic flow $\psi$ generated by a constant vector field, the group structure of $\hbox{\rm id}_{T_n}\to\hbox{\rm ker}\rho_\phi\to\rho_\phi^{-1}(\Lambda)\to\Lambda\to 1$ is determined by that of $\hbox{\rm id}_{T^n}\to\hbox{\rm ker}\rho_\psi\to\rho_\psi^{-1}(\Lambda)\to\Lambda\to 1$. Attention is therefore restricted to a quasiperiodic flow $\phi$ generated by a constant vector field $X$.

\section{Lifting the Generalized Symmetry Equation}

The equation $R_*X=\alpha X$ appearing in part c) of Theorem \ref{char} is the generalized symmetry equation for $\phi$. It is an equation on ${\bf T}T^n$, and lifting it to ${\bf T}{\mathbb R}^n$, the universal cover of ${\bf T}T^n$, requires lifting the diffeomorphism $R$ of $T^n$ to a diffeomorphism of ${\mathbb R}^n$, and the vector field $X$ on $T^n$ to a vector field on ${\mathbb R}^n$. The covering map $\pi:{\mathbb R}^n\to T^n$ is a local diffeomorphism for which
\[\pi(x+m)=\pi(x)\]
for any $x\in{\mathbb R}^n$ and any $m\in{\mathbb Z}^n$. Let $R:T^n\to T^n$ be a continuous map. A {\it lift} of $R\pi:{\mathbb R}^n\to T^n$ is a continuous map $Q:{\mathbb R}^n\to{\mathbb R}^n$ for which $R\pi=\pi Q$. Since $\pi$ is a fixed map, $Q$ is also said to be a lift of $R$. Any two lifts of $R$ differ by a deck transformation of $\pi$, which is a translation of ${\mathbb R}^n$ by an $m\in{\mathbb Z}^n$.

\begin{theorem}\label{lift3}
Let $R:T^n\to T^n$ and $Q:{\mathbb R}^n\to{\mathbb R}^n$. Then $Q$ is a lift of a diffeomorphism $R$ of $T^n$ if and only if $Q$ is a diffeomorphism of ${\mathbb R}^n$ such that

{\rm a)} for any $m\in{\mathbb Z}^n$, $Q(x+m)-Q(x)$ is independent of $x\in{\mathbb R}^n$, and

{\rm b)} the map $l_Q(m)=Q(x+m)-Q(x)$ is an isomorphism of ${\mathbb Z}^n$.
\end{theorem}

\begin{proof}
The proof uses standard arguments in topology.
\end{proof}

The canonical projections $\tau_{{\mathbb R}^n}:{\bf T}{\mathbb R}^n\to{\mathbb R}^n$ and $\tau_{T^n}:{\bf T}T^n\to T^n$ are smooth. The former is a lift of the latter,
\[\tau_{T^n}{\bf T}\pi=\pi\tau_{{\mathbb R}^n},\]
which lift sends $w\in{\bf T}_x{\mathbb R}^n$ to $x\in{\mathbb R}^n$. The covering map ${\bf T}\pi:{\bf T}{\mathbb R}^n \to {\bf T}T^n$ is a local diffeomorphism.

A vector field on $T^n$ is a smooth map $Y:T^n\to{\bf T}T^n$ such that $\tau_{T^n}Y=\hbox{\rm id}_{T^n}$. A vector field on ${\mathbb R}^n$ is a smooth map $Z:{\mathbb R}^n\to{\bf T}{\mathbb R}^n$ such that $\tau_{{\mathbb R}^n}Z=\hbox{\rm id}_{{\mathbb R}^n}$.

\begin{lemma}\label{lift4}
If $Y$ is a vector field on $T^n$, then there is only one lift of $Y$ that is a vector field on ${\mathbb R}^n$.
\end{lemma}

\begin{proof}
Let $x_0\in{\mathbb R}^n$, $\theta_0\in T^n$ be such that $Y\pi(x_0)=Y(\theta_0)$. Let $w_{x_0}\in{\bf T}_{x_0}{\mathbb R}^n$ be the only vector such that ${\bf T}\pi(w_{x_0})=Y(\theta_0)$. There exists a unique lift $Z:{\mathbb R}^n\to{\bf T}{\mathbb R}^n$ such that $Y\pi={\bf T}\pi Z$ and $Z(x_0)=w_{x_0}$ (Theorem 4.1, p. 143 \cite{BR}). Because $Y$ is a vector field on $T^n$, $Z$ is a lift of $Y\pi$, and $\tau_{{\mathbb R}^n}$ is a lift of $\tau_{T^n}$,
\[\pi(x)=\tau_{T_n}Y\pi(x)=\tau_{T^n}{\bf T}\pi Z(x)=\pi\tau_{{\mathbb R}^n}Z(x).\]
So the difference $x-\tau_{{\mathbb R}^n}Z(x)$ is a discrete valued map. Because ${\mathbb R}^n$ is connected, this difference is a constant (see Proposition 4.5, p. 10 \cite{BR}). This constant is zero because $\tau_{{\mathbb R}^n}Z(x_0)=x_0$, and so $\tau_{{\mathbb R}^n}Z=\hbox{\rm id}_{{\mathbb R}^n}$. The equation $Y\pi={\bf T}\pi Z$ implies that $Z$ is smooth because $\pi$ and ${\bf T}\pi$ are local diffeomorphisms and because $Y$ is smooth. The choice of the only vector $w\in{\bf T}_{x_0+m}{\mathbb R}^n$ for any $0\ne m\in{\mathbb Z}^n$ such that ${\bf T}\pi(w)=Y(\theta_0)$ would lead to a lift $Z_m$ of $Y$ that is not a vector field on ${\mathbb R}^n$ because $\tau_{{\mathbb R}^n}Z_m(x)=x+m$. The collection $\{Z_m:m\in{\mathbb Z}\}$, with $Z_0=Z$, accounts for all the lifts of $Y$ by the uniqueness of the lift and the uniqueness of the vector $w$. Therefore $Z$ is the only lift of $Y$ that is a vector field on ${\mathbb R}^n$.
\end{proof}

For a vector field $X$ on $T^n$, let $\hat X$ denote the only lift of $X$ that is a vector field on ${\mathbb R}^n$ as described in Lemma \ref{lift4}; $\hat X$ satisfies $X\pi={\bf T}\pi \hat X$. For a diffeomorphism $R$ of $T^n$, let $\hat R$ be a lift of $R$; the lift $\hat R$ is a diffeomorphism of ${\mathbb R}^n$ (by Theorem \ref{lift3}) for which $R\pi=\pi\hat R$.

\begin{lemma}\label{lift5}
The only lift of the vector field $R_*X$ on $T^n$ that is a vector field on ${\mathbb R}^n$ is $\hat R_*\hat X$.
\end{lemma}

\begin{proof}
A lift of $R_*X$ is $\hat R_*\hat X$ because
\begin{eqnarray*}
{\bf T}\pi\hat R_*\hat X
& = & {\bf T}\pi{\bf T}\hat R\hat X\hat R^{-1} \\
& = & {\bf T}(\pi \hat R)\hat X\hat R^{-1} \\
& = & {\bf T}(R\pi)\hat X\hat R^{-1} \\
& = & {\bf T}R{\bf T}\pi\hat X\hat R^{-1} \\
& = & {\bf T}RX\pi\hat R^{-1} \\
& = & {\bf T}RXR^{-1}\pi \\
& = & R_*X\pi.
\end{eqnarray*}
By definition, $\hat R_*\hat X$ is a vector field on ${\mathbb R}^n$. By Lemma \ref{lift4}, it is the only lift of $R_*X$ that is a vector field on ${\mathbb R}^n$.
\end{proof}

\begin{lemma}\label{lift6}
For any $\alpha\in{\mathbb R}^*$, the only lift of the vector field $\alpha X$ on $T^n$ that is a vector field on ${\mathbb R}^n$ is $\alpha\hat X$.
\end{lemma}

\begin{proof}
A lift of $\alpha X$ is $\alpha\hat X$ because 
\[ {\bf T}\pi(\alpha\hat X)=\alpha{\bf T}\pi\hat X=\alpha X\pi.\]
Only one lift of $\alpha X$ is a vector field (Lemma \ref{lift4}), and $\alpha\hat X$ is this lift.
\end{proof}

\begin{theorem}\label{lift7}
Let $X$ be a vector field on $T^n$, $\hat X$ the lift of $X$ that is a vector field on ${\mathbb R}^n$, $R$ a diffeomorphism of $T^n$, $\hat R$ a lift of $R$, and $\alpha$ a nonzero real number. Then $R_*X=\alpha X$ if and only if $\hat R_*\hat X=\alpha\hat X$.
\end{theorem}

\begin{proof}
Suppose that $R_*X=\alpha X$. By Lemma \ref{lift5}, $\hat R_*\hat X$ is a lift of $R_*X$: ${\bf T}\pi\hat R_*\hat X=R_*X\pi$. By Lemma \ref{lift6}, $\alpha \hat X$ is a lift of $\alpha X$: ${\bf T}\pi(\alpha \hat X)=\alpha X\pi$. Then
\[{\bf T}\pi\big(\hat R_*\hat X-\alpha\hat X\big)=\big(R_*X-\alpha X\big)\pi={\bf 0}_{T^n}\pi,\]
where ${\bf 0}_{T^n}$ is the zero vector field on $T^n$. So $\hat R_*\hat X-\alpha \hat X$ is a lift of ${\bf 0}_{T^n}$. The only lift of ${\bf 0}_{T^n}$ that is a vector field on ${\mathbb R}^n$ is ${\bf 0}_{{\mathbb R}^n}$, the zero vector field on ${\mathbb R}^n$. By Lemma \ref{lift5} and Lemma \ref{lift6}, the difference $\hat R_*\hat X-\alpha\hat X$ is a vector field on ${\mathbb R}^n$. By Lemma \ref{lift4}, $\hat R_*\hat X-\alpha\hat X={\bf 0}_{{\mathbb R}^n}$. Thus, $\hat R_*\hat X=\alpha\hat X$.

Suppose that $\hat R_*\hat X=\alpha\hat X$. Then
\begin{eqnarray*}
R_*X\pi
& = & {\bf T}RXR^{-1}\pi \\
& = & {\bf T}RX\pi\hat R^{-1} \\
& = & {\bf T}R{\bf T}\pi\hat X\hat R^{-1} \\
& = & {\bf T}(R\pi)\hat X\hat R^{-1} \\
& = & {\bf T}(\pi\hat R)\hat X\hat R^{-1} \\
& = & {\bf T}\pi{\bf T}\hat R\hat X\hat R^{-1} \\
& = & {\bf T}\pi\hat R_*\hat X \\
& = & {\bf T}\pi(\alpha\hat X) \\
& = & \alpha{\bf T}\pi\hat X \\
& = & \alpha X\pi.
\end{eqnarray*}
Because $\pi$ is surjective, $R_*X=\alpha X$.
\end{proof}

\section{Solving the Lifted Generalized Symmetry Equation}

The lift of $R_*X=\alpha X$ is an equation on ${\bf T}{\mathbb R}^n$ of the form $Q_*\hat X=\alpha\hat X$ for $Q\in\hbox{\rm Diff}({\mathbb R}^n)$. With global coordinates $x=(x_1,x_2,...,x_n)$ on ${\mathbb R}^n$, the diffeomorphism $Q$ has the form
\[Q(x_1,x_2,...,x_n)=(f_1(x_1,x_2,...,x_n),...,f_n(x_1,x_2,...,x_n))\]
for smooth functions $f_i:{\mathbb R}^n\to{\mathbb R}$, $i=1,...,n$. Let $\theta=(\theta_1,\theta_2,...,\theta_n)$ be global coordinates on $T^n$ such that $\theta_i=x_i$ mod 1, $i=1,2,...,n$. If
\[X(\theta)=a_1\frac{\partial}{\partial \theta_1}+ a_2\frac{\partial}{\partial \theta_2}+\cdot\cdot\cdot+a_n\frac{\partial}{\partial \theta_n}\]
for constants $a_i\in{\mathbb R}$, $i=1,...,n$, then
\[\hat X(x)= a_1\frac{\partial}{\partial x_1}+ a_2\frac{\partial}{\partial x_2}+\cdot\cdot\cdot+a_n\frac{\partial}{\partial x_n},\]
so that $Q_*\hat X=\alpha\hat X$ has the form
\[\sum^n_{j=1}a_j\frac{\partial f_i}{\partial x_j}=\alpha a_i,\ i=1,...,n.\]
This is an uncoupled system of linear, first order equations which is readily solved for its general solution.

\begin{lemma}\label{solve1}
For real numbers $a_1,a_2,...,a_n$ and $\alpha$ with $a_n\ne0$, the general solution of the system of $n$ linear partial differential equations
\[\sum^n_{j=1}a_j\frac{\partial f_i}{\partial x_j}=\alpha a_i,\ i=1,...,n\]
is
\[
f_i(x)=\alpha\frac{a_i}{a_n}x_n + h_i\bigg(x_1-\frac{a_1}{a_n}x_n,x_2-\frac{a_2}{a_n}x_n,...,x_{n-1}-\frac{a_{n-1}}{a_n}x_n\bigg),
\]
for arbitrary smooth functions $h_i:{\mathbb R}^{n-1}\to{\mathbb R}$, $i=1,...,n$.
\end{lemma}

\begin{proof}
For each $i=1,...,n$, consider the initial value problem
\begin{eqnarray*}
\sum^n_{j=1}a_j\frac{\partial f_i}{\partial x_j}&=&\alpha a_i \\
x_j(0,s_1,s_2,...,s_{n-1}) &=& s_j \hbox{\rm\ for\ }j=1,...,n-1 \\
x_n(0,s_1,s_2,...,s_{n-1}) &=& 0 \\
f_i(0,s_1,s_2,...,s_{n-1}) &=& h_i(s_1,s_2,...,s_{n-1})
\end{eqnarray*}
for parameters $(s_1,s_2,...,s_{n-1})\in{\mathbb R}^{n-1}$ and initial data $h_i:{\mathbb R}^{n-1}\to{\mathbb R}$. Using the method of characteristics (see \cite{JO} for example), the solution of the initial value problem in parametric form is
\begin{eqnarray*}
x_j(t,s_1,s_2,...,s_{n-1}) & = & a_jt+s_j\hbox{\rm\ for\ }j=1,...,n-1 \\
x_n(t,s_1,s_2,...,s_{n-1}) & = & a_nt \\
f_i(t,s_1,s_2,...,s_{n-1}) & = & \alpha a_it +h_i(s_1,s_2,...,s_{n-1}).
\end{eqnarray*}
The coordinates $(x_1,x_2,...,x_n)$ and the parameters $(t,s_1,s_2,...,s_{n-1})$ are related by
\[ \begin{bmatrix} x_1 \\ x_2 \\ x_3 \\ \vdots \\ x_{n-1} \\ x_n \end{bmatrix} =
\begin{bmatrix}
a_1 & 1 & 0 & 0 & \ldots & 0 \\
a_2 & 0 & 1 & 0 & \ldots & 0 \\
a_3 & 0 & 0 & 1 & \ldots & 0 \\
\vdots & \vdots & \vdots & \vdots & \ddots & \vdots \\
a_{n-1} & 0 & 0 & 0 & \ldots & 1 \\
a_n & 0 & 0 & 0 & \ldots & 0
\end{bmatrix}
\begin{bmatrix} t \\ s_1 \\ s_2 \\ \vdots \\ s_{n-2} \\ s_{n-1} 
\end{bmatrix}
\]
The determinant of the  $n\times n$ matrix is $(-1)^na_n$, which is nonzero by hypothesis. Inverting the matrix equation gives
\[ \begin{bmatrix} t \\ s_1 \\ s_2 \\ \vdots \\ s_{n-2} \\ s_{n-1} \end{bmatrix} =
\begin{bmatrix}
0 & 0 & \ldots & 0 & 0 & 1/a_n \\
1 & 0 & \ldots & 0 & 0 & -a_1/a_n \\
0 & 1 & \ldots & 0 & 0 & -a_2/a_n \\
\vdots & \vdots & \ddots & \vdots & \vdots & \vdots \\
0 & 0 & \ldots & 1 & 0 & -a_{n-2}/a_n \\
0 & 0 & \ldots & 0 & 1 & -a_{n-1}/a_n 
\end{bmatrix}
\begin{bmatrix}
x_1 \\ x_2 \\ x_3 \\ \vdots \\ x_{n-1} \\ x_n
\end{bmatrix}
\]
Substitution of the expressions for $t$ and the $s_i$'s in terms of the $x_i$'s into
\[f_i(x_1,x_2,...,x_n)=\alpha a_it+h_i(s_1,s_2,...,s_{n-1})\]
gives the desired form of the general solution. 
\end{proof}

\begin{lemma}\label{solve2}
If $a_1,a_2,...,a_n$ are independent over ${\mathbb Q}$, then
\[J=\bigg\{\bigg(m_1-\frac{a_1}{a_n}m_n,.. .,m_{n-1}-\frac{a_{n-1}}{a_n}m_n\bigg): m_1,.. .,m_n\in{\mathbb Z}\bigg\}\]
is a dense subset of ${\mathbb R}^{n-1}$.
\end{lemma}

\begin{proof}
Suppose $a_1,a_2,...,a_n$ are independent over ${\mathbb Q}$. This implies that none of the $a_i$'s are zero. In particular, $a_n\ne0$. Consider the flow
\[\psi_t(\theta_1,...,\theta_{n-1},\theta_n)=(\theta_1-(a_1/a_n)t,...,\theta_{n-1}-(a_{n-1}/a_n)t,\theta_n-t)\]
on $T^n$ which is generated by the vector field
\[Y=-\frac{a_1}{a_n}\frac{\partial}{\partial \theta_1}-\frac{a_2}{a_n}\frac{\partial}{\partial \theta_2}-\cdot\cdot\cdot-\frac{a_{n-1}}{a_n}\frac{\partial}{\partial \theta_{n-1}}-\frac{\partial}{\partial \theta_n}.\]
The coefficients of $Y$ are independent over ${\mathbb Q}$ because $a_1,a_2,...,a_n$ are independent over ${\mathbb Q}$ and
\[m_1a_1+\cdot\cdot\cdot+m_na_n=0\Leftrightarrow -m_1\frac{a_1}{a_n}-\cdot\cdot\cdot-m_{n-1}\frac{a_{n-1}}{a_n}-m_n=0.\]
So the orbit of $\psi$ through any point $\theta_0\in T^n$, \[\gamma_\psi(\theta_0)=\{\psi_t(\theta_0):t\in{\mathbb R}\},\]
is dense in $T^n$ (Corollary 1, p. 287 \cite{AR}).

The submanifold
\[ P=\{(\theta_1,...,\theta_{n-1},\theta_n):\theta_n=0\}\]
of $T^n$, which is diffeomorphic to $T^{n-1}$, is a global Poincar\'e section for $\psi$ because $X(\theta)\not\in{\bf T}_\theta P$ for every $\theta\in P$ and because $\gamma_\psi(\theta_0)\cap P\ne\emptyset$ for every $\theta_0\in T^n$. Define the projection $\wp:T^n \to T^{n-1}$ by
\[\wp(\theta_1,\theta_2,...,,\theta_{n-1},\theta_n)=(\theta_1,\theta_2,...,\theta_{n-1})\]
and the injection $\imath:T^{n-1}\to T^n$ by
\[\imath(\theta_1,\theta_2,...,\theta_{n-1})=(\theta_1,\theta_2,...,\theta_{n-1},0).\]
The Poincar\'e map induced on $\wp(P)$ by $\psi$ is given by $\bar\psi=\wp\psi_1\imath$ because $\psi_1(\theta_0)\in P$ when $\theta_0\in P$. For any $\kappa\in{\mathbb Z}$, $\bar\psi^\kappa=\wp\psi_{\kappa}\imath$. So, for instance, with $0=(0,0,...,0)\in T^n$ and $\bar 0=\wp(0)$,
\begin{eqnarray*}
\wp\big(\gamma_\psi(0)\cap P\big)
& = & \{\bar\psi^\kappa(\bar 0):\kappa\in{\mathbb Z}\} \\
& = & \bigg\{\bigg(-\frac{a_1}{a_n}\kappa,-\frac{a_2}{a_n}\kappa,\cdot\cdot\cdot,-\frac{a_{n-1}}{a_n}\kappa\bigg):\kappa\in{\mathbb Z}\bigg\},
\end{eqnarray*}
where for each $i=1,...,n-1$, the quantity $-(a_i/a_n)\kappa$ is taken mod 1. With $\bar\pi:{\mathbb R}^{n-1}\to T^{n-1}$ as the covering map,
\[J=\bar\pi^{-1}\big(\wp(\gamma_\psi(0)\cap P)\big).\]
If $\wp(\gamma_\psi(0)\cap P)$ were dense in $\wp(P)$, then $J$ would be dense in $R^{n-1}$ because $\bar\pi$ is a covering map. (That is, if $\wp(\gamma_\psi(0)\cap P)\cap[0,1)^{n-1}$ is dense in the fundamental domain $[0,1)^{n-1}$ of the covering map $\bar\pi$, then by translation, it is dense in ${\mathbb R}^{n-1}$.)

Define $\chi:{\mathbb R}\times T^{n-1}\to T^n$ by
\[\chi(t,\theta_1,\theta_2,...,\theta_{n-1})=\psi\big(t,\imath(\theta_1,\theta_2,...,\theta_{n-1})\big).\]
The map $\chi$ is a local diffeomorphism by the Inverse Function Theorem because
\[{\bf T}\chi=
\begin{bmatrix}
-a_1/a_n & 1 & 0 & \ldots & 0 \\
-a_2/a_n & 0 & 1 & \ldots & 0 \\
\vdots & \vdots & \vdots & \ddots & \vdots \\
-a_{n-1}/a_n & 0 & 0 &\ldots & 1 \\
-1 & 0 & 0 & \ldots & 0
\end{bmatrix}\]
has determinant of $(-1)^{n+1}$. Let $O$ be a small open subset of $\wp(P)$. For $\epsilon>0$, the set $O_\epsilon=(-\epsilon,\epsilon)\times O$ is an open subset in the domain of $\chi$. For $\epsilon$ small enough, the image $\chi(O_\epsilon)$ is open in $T^n$ because $\chi$ is a local diffeomorphism. By the denseness of $\gamma_\psi(0)$ in $T^n$, there is a point $\theta_0$ in $\chi(O_\epsilon)\cap\gamma_\psi(0)$. By the definition of $\chi(O_\epsilon)$, there is an $\bar\epsilon\in(-\epsilon,\epsilon)$ and a $\bar\theta_0\in O$ such that $\chi(\bar\epsilon,\bar\theta_0)=\theta_0$. Thus $\imath(\bar \theta_0)\in\gamma_\psi(0)$, and so $\wp(\gamma_\psi(0)\cap P)$ intersects $O$ at $\bar\theta_0$. Since $O$ is any small open subset of $\wp(P)$, the set $\wp(\gamma_\psi(0)\cap P)$ is dense in $\wp(P)$.
\end{proof}

\begin{theorem}\label{solve3}
If $\alpha\in{\mathbb R}^*$ and the coefficients of $X=\sum_{i=1}^n a_i\partial/\partial\theta_i$ are independent over ${\mathbb Q}$, then for each $R\in\hbox{\rm Diff}(T^n)$ that satisfies $R_*X=\alpha X$ there exist $B=(b_{ij})\in\hbox{\rm GL}(n,{\mathbb Z})$ and $c\in{\mathbb R}^n$ such that
\[\hat R(x)=Bx+c\]
for $x=(x_1,x_2,...,x_n)$, in which
\[b_{in}= \alpha\frac{a_i}{a_n}-\sum_{j=1}^{n-1}b_{ij}\frac{a_j}{a_n},\ i=1,...,n.\]
\end{theorem}

\begin{proof}
Suppose that the $a_1,a_2,...,a_n$ are independent over ${\mathbb Q}$. For $\alpha\in{\mathbb R}^*$, suppose that $R\in\hbox{\rm Diff}(T^n)$ is a solution of $R_*X=\alpha X$. A lift $\hat R$ of $R$ is a diffeomorphism of ${\mathbb R}^n$ by Theorem \ref{lift3}. The lift of $X$ that is a vector field on ${\mathbb R}^n$ is $\hat X=\sum^n_{i=1}a_i(\partial/\partial x_i)$. By Theorem \ref{lift7}, $\hat R$ is a solution of $\hat R_*\hat X=\alpha\hat X$. With global coordinates $(x_1,x_2,...,x_n)$ on ${\mathbb R }^n$ write 
\[\hat R(x)=(f_1(x_1,...,x_n),...,f_n(x_1,...,x_n)).\]
In terms of this coordinate description, the equation $\hat R_*\hat X=\alpha X$ written out is
\[\sum_{j=1}^n a_j\frac{\partial f_i}{\partial x_j}=\alpha a_i,\ i=1,...,n.\]
The independence of the coefficients of $\hat X$ over ${\mathbb Q}$ implies that $a_n\ne0$. By Lemma \ref{solve1}, there are smooth functions $h_i:{\mathbb R}^{n-1}\to{\mathbb R}$, $i=1,...,n$, such that
\[
f_i(x_1,...,x_n)=\alpha\frac{a_i}{a_n}x_n+h_i(s_1,s_2,...,s_{n-1})\]
where
\[s_i=x_i-\frac{a_i}{a_n}x_n,\ i=1,...,n-1.\]

By Theorem \ref{lift3}, $\hat R(x+m)-\hat R(x)$ is independent of $x$ for each $m\in{\mathbb R}^n$. This implies for each $i=1,...,n$ that
\begin{eqnarray*}
& & f_i(x+m)-f_i(x)\\
& & \ \ =f_i(x_1+m_1,x_2+m_2,...,x_n+m_n)-f_i(x_1,x_2,...,x_n) \\
& & \ \ =\alpha\frac{a_i}{a_n}m_n + h_i\bigg(s_1+m_1-\frac{a_1}{a_n}m_n,...,s_{n-1}+m_{n-1}-\frac{a_{n-1}}{a_n}m_n\bigg)\\
& & \ \ \ \ - h_i(s_1,...,s_{n-1})
\end{eqnarray*}
is independent of $x$ for every $m=(m_1,m_2,...,m_n)\in{\mathbb Z}^n$. This independence means that $f_i(x+m)-f_i(x)$ is a function of $m$ only. So for each $j=1,.. .,n-1$,
\begin{eqnarray*}
0 & = & \frac{\partial}{\partial x_j}\big[f_i(x_1+m_1,x_2+m_2,. ..,x_n+m_n)-f_i(x_1,x_2,...,x_n)\big] \\
& = & \frac{\partial h_i}{\partial s_j}\bigg(s_1+m_1-\frac{a_1}{a_n}m_n,...,s_{n-1}+m_{n-1}-\frac{a_{n-1}}{a_n}m_n\bigg) \\
&  & -\frac{\partial h_i}{\partial s_j}\big(s_1,...,s_{n-1}\big).
\end{eqnarray*}
So, in particular
\[\frac{\partial h_i}{\partial s_j}\bigg(m_1-\frac{a_1}{a_n}m_n,...,m_{n-1}-\frac{a_{n-1}}{a_n}m_n\bigg)= 
\frac{\partial h_i}{\partial s_j}\big(0,...,0\big) \]
for all $(m_1,...,m_n)\in{\mathbb Z}^n$. By Lemma \ref{solve2}, the set
\[\bigg\{\bigg(m_1-\frac{a_1}{a_n}m_n,.. .,m_{n-1}-\frac{a_{n-1}}{a_n}m_n\bigg): m_1,.. .,m_n\in{\mathbb Z}\bigg\}\]
is dense in ${\mathbb R}^{n-1}$, which together with the smoothness of $h_i$ implies that $\partial h_i/\partial s_j$ is a constant. Let this constant be $b_{ij}$ for $i=1,...,n$, $j=1,...,n-1$. By Taylor's Theorem,
\[h_i(s_1,...,s_{n-1})=c_i + \sum_{j=1}^{n-1} b_{ij}s_j\]
for constants $c_i\in{\mathbb R}$. Thus,
\begin{eqnarray*}
f_i(x_1,...,x_n)
& = & c_i + \alpha\frac{a_i}{a_n}x_n+\sum_{j=1}^{n-1} b_{ij}\bigg(x_j-\frac{a_j}{a_n}x_n\bigg) \\
& = & c_i + \sum_{j=1}^{n-1}b_{ij}x_j+\bigg(\alpha\frac{a_i}{a_n}-\sum_{j=1}^{n-1}b_{ij}\frac{a_j}{a_n}\bigg)x_n.
\end{eqnarray*}
For each $i=1,2,...,n$, set
\[b_{in}=\alpha\frac{a_i}{a_n}-\sum_{j=1}^{n-1}b_{ij}\frac{a_j}{a_n}\bigg.\]
Then for each $i=1,2,...,n$,
\[f_i(x_1,x_2,...,x_n)=c_i + \sum_{j=1}^n b_{ij}x_j.\]
So $\hat R$ has the form $\hat R(x)=Bx+c$ where $B=(b_{ij})$ is an $n\times n$ matrix, and $c\in{\mathbb R}^n$.

By Theorem \ref{lift3}, the map $l_{\hat R}(m)=\hat R(x+m)-\hat R(x)$ is an isomorphism of ${\mathbb Z}^n$. By the formula for $f_i$ derived above,
\[f_i(x_1+m_1,...,x_n+m_n)-f_i(x_1,x_2,...,x_m)=\sum_{j=1}^n b_{ij}m_j\]
for each $i=1,2,...,n$. This implies that
\[l_{\hat R}(m)=Bm.\]
Since $l_{\hat R}$ is an isomorphism of ${\mathbb Z}^n$, it follows that $B\in\hbox{\rm GL}(n,{\mathbb Z})$.
\end{proof}

Theorem \ref{solve3} restricts the search for lifts of generalized symmetries of a quasiperiodic flow on $T^n$ to affine maps on ${\mathbb R}^n$ of the form $Q(x)=Bx+c$ for $B\in{\rm GL}(n,{\mathbb Z})$ and $c\in{\mathbb R}^n$. For an affine map of this form, the difference
\[Q(x+m)-Q(x)=B(x+m)+c-(Bx+c)=Bm\]
is independent of $x$, and the map $l_Q(m)=Q(x+m)-Q(x)$ is an isomorphism of ${\mathbb Z}^n$, so that $Q$ is a lift of a diffeomorphism $R$ on $T^n$ by Theorem \ref{lift3}. If $Q$ is a solution of $Q_*\hat X=\alpha\hat X$, then by Theorem \ref{lift7}, $R$ is a solution of $R_*X=\alpha X$, so that by Theorem \ref{char}, $R\in S_\phi$.

Theorem \ref{solve3} also restricts the possibilities for the multipliers of any generalized symmetries of a quasiperiodic flow on $T^n$. One restriction employs the notion of an {\it algebraic integer}, which is a complex number that is a root of a monic polynomial in the polynomial ring ${\mathbb Z}[z]$. If $m$ is the smallest degree of a monic polynomial in ${\mathbb Z}[z]$ for which an algebraic integer is a root, then $m$ is the {\it degree} of that algebraic number (Definition 1.1, p. 1 \cite{RM}).

\begin{corollary}\label{solve4}
If the coefficients of $X=\sum_{i=1}^n a_i \partial/\partial\theta_i$ are independent over ${\mathbb Q}$, then each $\alpha\in\rho_\phi(S_\phi)$ is a real algebraic integer of degree at most $n$, and $\rho_\phi(S_\phi)\cap{\mathbb Q}=\{1,-1\}.$
\end{corollary}

\begin{proof}
For each $\alpha\in\rho_\phi(S_\phi)$ (which is real) there is an $R\in S_\phi$ such that $\rho_\phi(R)=\alpha$. By Theorem \ref{solve3} there is a $B\in\hbox{\rm GL}(n,{\mathbb Z})$ such that ${\bf T}\hat R=B$. Then by Theorem \ref{char} and Theorem \ref{lift7},
\[B\hat X=\hat R_*\hat X=\alpha\hat X.\]
So, $\alpha$ is an eigenvalue of $B$ (and $\hat X$ is an eigenvector of $B$.) The characteristic polynomial of $B$ is an n-degree monic polynomial in ${\mathbb Z}[z]$:
\[z^n+d_{n-1}z^{n-1}+\cdot\cdot\cdot+d_1z+d_0.\]
Thus $\alpha$ is a real algebraic integer of degree at most $n$. The value of $d_0$ is $\hbox{\rm det}(B)$, which is a unit in ${\mathbb Z}$ (Theorem 3.5, p. 351 \cite{HU}). The only units in ${\mathbb Z}$ are $\pm 1$. So the only possible rational roots of the characteristic polynomial of $B$ are $\pm1$ (Proposition 6.8, p. 160 \cite{HU}). This means that $\rho_\phi(S_\phi)\cap{\mathbb Q}\subset\{1,-1\}$. But $\rho_\phi(S_\phi)\cap{\mathbb Q}\supset\{1,-1\}$ by Theorem \ref{mult}. Thus, $\rho_\phi(S_\phi)\cap{\mathbb Q}=\{1,-1\}$.
\end{proof}

Another restriction on the possibilities for the multipliers of any generalized symmetries of $\phi$ employs linear combinations over ${\mathbb Z}$ of pair wise ratios of the entries of the ``eigenvector'' $\hat X$.

\begin{corollary}\label{solve5}
If the coefficients of $X=\sum^n_{i=1}a_i\partial/\partial\theta_i$ are independent over ${\mathbb Q}$, then for any $\alpha\in\rho_\phi(S_\phi)$ there exists a $B=(b_{ij})\in\hbox{\rm GL}(n,{\mathbb Z})$ such that
\[\alpha=\sum_{j=1}^n b_{ij}\frac{a_j}{a_i}\]
for each $i=1,...,n$.
\end{corollary}

\begin{proof}
Suppose that $\alpha\in\rho_\phi(S_\phi)$. Then there is an $R\in S_\phi$ such that $\alpha=\rho_\phi(R)$. By Theorem \ref{solve3}, there is a $B=(b_{ij})\in{\rm GL}(n,{\mathbb Z})$ such that ${\bf T}\hat R=B$ with
\[b_{in}= \alpha\frac{a_i}{a_n}-\sum_{j=1}^{n-1}b_{ij}\frac{a_j}{a_n}\]
for each $i=1,...,n$. Solving this equation for $\alpha$ gives
\[\alpha=\sum_{j=1}^n b_{ij}\frac{a_j}{a_i}\]
for each $i=1,...,n$.
\end{proof}

The multiplier group of any quasiperiodic flow $\phi$ always contains $\{1,-1\}$ as stated in Theorem \ref{mult}. For each $t\in{\mathbb R}$, the diffeomorphism $\phi_t$ is in $S_\phi$ by definition. A lift of $\phi_t$ is $\hat \phi_t(x)=Ix+t\hat X$, where $I=\delta_{ij}$ is the $n\times n$ identity matrix, so that by Corollary \ref{solve5},
\[\alpha=\sum^n_{j=1}\delta_{ij}\frac{a_j}{a_i}=\frac{a_i}{a_i}=1\]
for each $i=1,...,n$. A lift of the reversing involution $N$ defined in the proof of Theorem \ref{mult} is $\hat N(x)=-Ix$, so that by Corollary \ref{solve5},
\[\alpha=-\sum_{j=1}^n \delta_{ij}\frac{a_j}{a_i}=-\frac{a_i}{a_i}=-1\]
for each $i=1,...,n$. Corollary \ref{solve5} enables a complete description of all symmetries and reversing symmetries of $\phi$.

\begin{theorem}\label{solve6}
Suppose that the coefficients of $X=\sum^n_{i=1}a_i\partial/\partial\theta_i$ are independent over ${\mathbb Q}$. If $\rho_\phi(R)=\pm1$ for an $R\in S_\phi$, then there is a $c\in{\mathbb R}^n$ such that $\hat R(x)=\rho_\phi(R)Ix+c$.
\end{theorem}

\begin{proof}
Let $R\in S_\phi$. By Theorem \ref{solve3} there exists a $B=(b_{ij})\in\hbox{\rm GL}(n,{\mathbb Z})$ and a $c\in{\mathbb R}^n$ such that $\hat R(x)=Bx+c$. By Corollary \ref{solve5}, the entries of $B$ satisfy
\[\rho_\phi(R)=\sum_{j=1}^n b_{ij}\frac{a_j}{a_i}\]
for each $i=1,2,...,n$. By hypothesis, $\rho_\phi(R)=\pm 1$. Then for each $i=1,2,...,n$,
\[b_{i1}a_1+\cdot\cdot\cdot+(b_{ii}\mp1)a_i+\cdot\cdot\cdot+b_{in}a_n=0.\]
By the independence of $a_1,a_2,...,a_n$ over ${\mathbb Q}$, $b_{ij}=0$ when $i\ne j$ and $b_{ii}=\rho_\phi(R)$ for all $i=1,2,...,n$. Therefore, $\hat R(x)=\rho_\phi(R)Ix+c$.
\end{proof}

\begin{corollary}\label{solve7}
If the coefficients of $X=\sum^n_{i=1}a_i\partial/\partial\theta_i$ are independent over ${\mathbb Q}$, then $\hbox{\rm ker}\rho_\phi\cong T^n$.
\end{corollary}

\begin{proof}
Let $R\in S_\phi$ such that $\rho_\phi(R)=1$. By Theorem \ref{solve6}, $\hat R(x)=Ix+c$ for some $c\in{\mathbb R}^n$. Now, for any $c\in{\mathbb R}^n$, the $Q\in\hbox{\rm Diff}(T^n)$ induced by $\hat Q(x)=Ix+c$ satisfies $Q_*X=X$ by Theorem \ref{lift7} because $\hat Q_*\hat X=\hat X$. So, by Theorem \ref{char}, $Q\in\hbox{\rm ker}\rho_\phi$. Since $c$ is arbitrary, $Q\pi=\pi\hat Q$, and $\pi({\mathbb R}^n)=T^n$, it follows that $\hbox{\rm ker}\rho_\phi\cong T^n$. 
\end{proof}

\begin{corollary}\label{solve8}
If the coefficients of $X=\sum^n_{i=1}a_i\partial/\partial\theta_i$ are independent over ${\mathbb Q}$, then every reversing symmetry of $\phi$ is a reversing involution.
\end{corollary}

\begin{proof}
Suppose $R\in S_\phi$ is a reversing symmetry. By Theorem \ref{solve6}, $\hat R(x)=-Ix+c$ for some $c\in{\mathbb R}^n$, and so $\hat R^2(x)=Ix$. This implies that $R^2=\hbox{\rm id}_{T^n}$.
\end{proof}

\section{A Splitting Map for the Extension}

A splitting map for the short exact sequence,
\[ \hbox{\rm id}_{T^n}\to\hbox{\rm ker}\rho_\phi\to \rho_\phi^{-1}(\Lambda)\to\Lambda\to 1,\]
is a homomorphism $h_\Lambda:\Lambda\to\rho_\phi^{-1}(\Lambda)$ such that $j_\Lambda h_\Lambda=1$, the identity isomorphism of $\Lambda$. Such a splitting map exists if and only if $\rho_\phi^{-1}(\Lambda)$ is the internal semidirect product of $\hbox{\rm ker}\rho_\phi$ by a subgroup isomorphic to $\Lambda$ (Theorem 9.5.1, p. 240 \cite{SC}).

Given that $\phi$ is generated by a constant vector field, set
\[ \Pi_\phi=\{B\in\hbox{\rm GL}(n,{\mathbb Z}):B={\bf T}\hat R\hbox{\rm\ for an\ }R\in S_\phi\},\]
and define a map $\nu_\phi:\Pi_\phi\to\rho_\phi(S_\phi)$ by $\nu_\phi(B)=\rho_\phi(R)$ where $R\in S_\phi$ with ${\bf T}\hat R=B$.

\begin{lemma}\label{rep0}
If the coefficients of $X=\sum^n_{i=1}a_i\partial/\partial\theta_i$ are independent over ${\mathbb Q}$, then $\nu_\phi$ is well-defined.
\end{lemma}

\begin{proof}
Let $B\in\Pi_\phi$, and suppose there are $R,Q\in S_\phi$ with ${\bf T}\hat R=B={\bf T}\hat Q$ such that $\nu_\phi(B)=\rho_\phi(R)$ and $\nu_\phi(B)=\rho_\phi(Q)$. Then $RQ^{-1}\in S_\phi$ and $\hat R\hat Q^{-1}$ is a lift of $RQ^{-1}$ for which ${\bf T}(\hat R\hat Q^{-1})=BB^{-1}=I$. Hence $\hat R\hat Q^{-1}(x)=Ix+c$ for some $c\in{\mathbb R}^n$. This implies that $(\hat R\hat Q^{-1})_*\hat X=\hat X$, so that by Theorem \ref{lift7}, $(RQ^{-1})_*X=X$. By Theorem \ref{char}, $\rho_\phi(RQ^{-1})=1$. Because $\rho_\phi$ is a homomorphism, $\rho_\phi(R)=\rho_\phi(Q)$.
\end{proof}

\begin{lemma}\label{rep1}
If the coefficients of $X=\sum^n_{i=1}a_i\partial/\partial\theta_i$ are independent over ${\mathbb Q}$, then $\Pi_\phi$ is a subgroup of $\hbox{\rm GL}(n,{\mathbb Z})$.
\end{lemma}

\begin{proof}
Let $B,C\in\Pi_\phi$. Then there are $R,Q\in S_\phi$ such that ${\bf T}\hat R=B$ and ${\bf T}\hat Q=C$. The latter implies that ${\bf T}\hat Q^{-1}=({\bf T}\hat Q)^{-1}=C^{-1}$. Then $BC^{-1}={\bf T}\hat R{\bf T}\hat Q^{-1}={\bf T}(\hat R\hat Q^{-1})$. The diffeomorphism $x\to \hat R\hat Q^{-1}x$ of ${\mathbb R}^n$ satisfies conditions a) and b) of Theorem \ref{lift3}, and so is a lift of a diffeomorphism $V$ of $T^n$. Let $\alpha=\rho_\phi(R)$ and $\beta=\rho_\phi(Q)$. Then $\rho_\phi(Q^{-1})=\beta^{-1}$ because $\rho_\phi$ is a homomorphism, and so $(\hat Q^{-1})_*\hat X=\beta^{-1}\hat X$. Thus,
\[{\bf T}(\hat R\hat Q^{-1})\hat X=(\hat R\hat Q^{-1})_*\hat X=\alpha\beta^{-1}\hat X.\]
By Theorem \ref{lift7}, $V_*X=\alpha\beta^{-1}X$, so that by Theorem \ref{char}, $V\in S_\phi$. The lifts $\hat R\hat Q^{-1}$ and $\hat V$ of $V$ differ by a deck transformation of $\pi$, so that $BC^{-1}={\bf T}(\hat R\hat Q^{-1})={\bf T}\hat V$. Therefore, $BC^{-1}\in \Pi_\phi$.
\end{proof}

\begin{theorem}\label{rep2}
If the coefficients of $X=\sum^n_{i=1}a_i\partial/\partial\theta_i$ are independent over ${\mathbb Q}$, then $\nu_\phi$ is an isomorphism and $\Pi_\phi$ is an abelian subgroup of $\hbox{\rm GL}(n,{\mathbb Z})$.
\end{theorem}

\begin{proof}
Let $B,C\in\Pi_\phi$. Then there are $R,Q\in S_\phi$ such that ${\bf T}\hat R=B$ and ${\bf T}\hat Q=C$. Let $\alpha=\rho_\phi(R)$ and $\beta=\rho_\phi(Q)$. By Theorem \ref{char} and Theorem \ref{lift7}, ${\bf T}\hat R\hat X=\alpha\hat X$ and ${\bf T}\hat Q\hat X=\beta\hat X$. By Lemma \ref{rep1}, $BC\in \Pi_\phi$, so that there is a $V\in S_\phi$ such that ${\bf T}\hat V=BC$. Hence,
\[\hat V_*\hat X={\bf T}\hat V\hat X=BC\hat X=\alpha\beta\hat X.\]
By Theorem \ref{lift7} and Theorem \ref{char}, $\rho_\phi(V)=\alpha\beta$. Thus, 
\[\nu_\phi(BC)=\alpha\beta=\nu_\phi(B)\nu_\phi(C).\]
By definition, $\nu_\phi$ is surjective, and by Theorem \ref{solve6}, $\hbox{\rm ker}\nu_\phi=\{I\}$. Therefore $\nu_\phi$ is an isomorphism. The multiplier group $\rho_\phi(S_\phi)$ is abelian because it is a subgroup of the abelian group ${\mathbb R}^*$. Thus $\Pi_\phi$ is abelian.
\end{proof}

The isomorphism $\nu_\phi:\Pi_\phi\to\Lambda$ enables the definition of a map $h_\Lambda:\Lambda\to\rho_\phi^{-1}(\Lambda)$. It is $h_\Lambda(\alpha)=R$ where $\hat R(x)=\nu_\phi^{-1}(\alpha)x$.

\begin{theorem}\label{rep4}
If the coefficients of $X=\sum^n_{i=1}a_i\partial/\partial\theta_i$ are independent over ${\mathbb Q}$, then $h_\Lambda$ is a splitting map for the extension $\hbox{\rm id}_{T^n}\to\hbox{\rm ker}\rho_\phi\to\rho_\phi^{-1}(\Lambda)\to\Lambda\to1$ for each $\Lambda<\rho_\phi(S_\phi)$.
\end{theorem}

\begin{proof}
For arbitrary $\alpha,\beta\in\Lambda$, set $R=h_\Lambda(\alpha)$, $Q=h_\Lambda(\beta)$, and $V=h_\Lambda(\alpha\beta)$. Then $\hat R(x)=\nu_\phi^{-1}(\alpha)x$, $\hat Q(x)=\nu_\phi^{-1}(\beta)x$, and $\hat V(x)=\nu_\phi^{-1}(\alpha\beta)x$. By Theorem \ref{rep2}, $\nu_\phi^{-1}$ is an isomorphism, so that $\hat V(x)=\nu_\phi^{-1}(\alpha)\nu_\phi^{-1}(\beta)x$. Because
\begin{eqnarray*}
h_\Lambda(\alpha)h_\Lambda(\beta)\pi(x)
& = & RQ\pi(x) \\
& = & \pi \hat R\hat Q(x) \\
& = & \pi\nu_\phi^{-1}(\alpha)\nu_\phi^{-1}(\beta)x \\
& = & \pi\nu_\phi^{-1}(\alpha\beta)x \\
& = & \pi \hat V(x) \\
& = & V\pi(x) \\
& = & h_\Lambda(\alpha\beta)\pi(x),
\end{eqnarray*}
and because $\pi$ is surjective, $h_\Lambda(\alpha)h_\Lambda(\beta)=h_\Lambda(\alpha\beta)$. Let $B={\bf T}\hat R=\nu_\phi^{-1}(\alpha)$. Then $\nu_\phi(B)=\rho_\phi(R)$, so that
\[ j_\Lambda h_\Lambda(\alpha)=j_\Lambda(R)=\rho_\phi(R)=\nu_\phi(B)=\nu_\phi(\nu_\phi^{-1}(\alpha))=\alpha.\]
Therefore $h_\Lambda$ is a splitting map for the extension.
\end{proof}

\begin{theorem}\label{rep5}
If the coefficients of $X=\sum^n_{i=1}a_i\partial/\partial\theta_i$ are independent over ${\mathbb Q}$, then for any $\{1\}\ne\Lambda<\rho_\phi(S_\phi)$, $\rho_\phi^{-1}(\Lambda)$ is the internal semidirect product of a subgroup isomorphic to the abelian group $T^n$ by a subgroup isomorphic to a multiplicative abelian group of real algebraic integers of degree at most $n$, which multiplicative group is isomorphic to a abelian subgroup of $\hbox{\rm GL}(n,{\mathbb Z})$.
\end{theorem}

\begin{proof}
By Theorem \ref{rep4}, $\rho_\phi^{-1}(\Lambda)$ is the internal semidirect product of $\hbox{\rm ker}\rho_\phi$ by $h_\Lambda(\Lambda)$. By Corollary \ref{solve7}, $\hbox{\rm ker}\rho_\phi\cong T^n$, and by Corollary \ref{solve4}, $\Lambda$ is a multilplicative group of algebraic integers of degree at most $n$, which by Theorem \ref{rep2} is isomorphic to an abelian subgroup of $\hbox{\rm GL}(n,{\mathbb Z})$. The internal semidirect product is not the internal direct product of $\hbox{\rm ker}\rho_\phi$ by $h_\Lambda(\Lambda)$ because the internal direct product of the two abelian groups is abelian while by Theorem \ref{mult}, $\rho_\phi^{-1}(\Lambda)$ is nonabelian whenever $\Lambda\ne\{1\}$.
\end{proof}

\end{document}